\documentclass[11pt]{article}

\usepackage{amsmath}
\usepackage{amsfonts}
\usepackage{amssymb}
\usepackage{amsthm}
\usepackage{graphics}
\usepackage{amscd}
\usepackage{graphicx,epsfig}
\usepackage{color}
\usepackage[all]{xy}
\usepackage{url}

\DeclareMathOperator{\Div}{div}

\renewcommand{\epsilon}{\varepsilon}

\newcommand{\boG}{\mathcal{G}}

\newcommand{\R}{\mathbb{R}}

\newcommand{\N}{\mathbb{N}}

\newcommand{\eps}{\varepsilon}

\newcommand{\Ome}{\Omega}

\newtheorem{thm}{Theorem}
\newtheorem{cor}[thm]{Corollary}
\newtheorem{prop}[thm]{Proposition}
\newtheorem{lem}[thm]{Lemma}
\newtheorem*{th1'}{Theorem 2'}
\newcommand{\inter}[1]{\overset{\circ}{#1}}

\renewcommand{\phi}{\varphi}

\newcounter{remark}

\newcounter{case}

\newcounter{construction}

\newcounter{fact}

\newenvironment{fact}[1]{\refstepcounter{fact}\label{#1}
  \noindent\textbf{Fact \ref{#1}.}}{}

\title{A height estimate for constant mean curvature graphs and uniqueness}
\author{Laurent Mazet}
\date{}

\begin{document}

\maketitle

\begin{abstract}
In this paper, we give a height estimate for constant mean curvature
graphs. Using this result we prove two results of uniqueness for the Dirichlet 
problem associated to the constant mean curvature equation on unbounded
domains. 
\end{abstract}


\section*{Introduction}
The surfaces with constant mean curvature are the mathematical modelling of
soap films. These surfaces appear as the interfaces in isoperimetric
problems. There exist different points of view on constant mean curvature
surfaces, one is to consider them as graphs.

Let $\Ome$ be a domain of $\R^2$. The graph of a function $u$ over $\Ome$ has
constant mean curvature $H>0$ if it satisfies the following partial 
differential equation:
\begin{equation*}\label{cmc}
\Div\left(\frac{\nabla u}{\sqrt{1+|\nabla u|^2}}\right)=2H
\tag{CMC}
\end{equation*}
The graph of such a solution is called a $H$-graph and has a upward pointing
mean curvature vector. 

The Dirichlet problem is a natural question about this point of view. For
bounded domains, after the work of J.~Serrin \cite{Se1}, J.~Spruck has given
in \cite{Sp} a general answer to the existence and uniqueness questions. His
results are of Jenkins-Serrin type \cite{JS} since infinite data are allowed. 

On unbounded domains, there is few constructions of solutions. The examples are
due to P.~Collin \cite{Co} and R.~Lop\'ez \cite{Lo} for graphs over a strip and
R.~Lop\'ez \cite{Lo,Lo2} for graphs with $0$ boundary data. 

In this paper, we investigate the uniqueness question of the Dirichlet
problem. To get uniqueness, we need a control of the solutions of the
Dirichlet problem, then we shall be able to bound the distance between two
solutions with the same boundary data.

To get this control, we use our main result which is Theorem \ref{estim}. We
call this result a height estimate since it bounds the difference of height
between two components of boundary of a $H$-graph. The idea of the proof is
that if the difference of height is too big, we can get a sphere of radius
$\frac{1}{H}$ through the $H$-graph; this is impossible because of the maximum
principle. 

With this height estimate, we can bound the difference between two solutions
of the same dirichlet problem under certain hypotheses on the boundary
value. For example, if we are on a strip, we prove that, if the boundary data
is Lipschitz, we have uniqueness of possible solutions. 

The first section of the paper is devoted to the statement of the height
estimate and its proof.

In Section \ref{conseq}, we give consequences of the height estimate for
solutions of the constant mean curvature equation on unbounded domains.

In the last section, we prove two uniqueness results for the Dirichlet problem
on unbounded domains.


\section{The height estimate}
\label{estimH1}
In this first section, we give a height estimate for solutions of the constant
mean curvature equation \eqref{cmc}. This estimate is designed for solutions
on unbounded domains.

First, we need some notations and remarks. If $\Ome$ is a domain in $\R^2$ and
$u$ is a function which is defined on $\Ome$, we note $F:\Ome\rightarrow \R^2$
the map with:
$$
F(x,y)=\big(x,u(x,y)\big)
$$

Let us explain what kind of domain we shall consider in the following.
For $a>\frac{1}{H}$ and $b>0$, we note $R_{a,b}=[-a,a]\times[-b,b]$. Let
$\Ome\subset R_{a,b}$ be a domain with piecewise smooth boundary. We
suppose that $\Ome$ satisfies the three following hypotheses :
\begin{enumerate}
\item $\Ome$ is connected, 
\item $\partial\Ome\cap\{-a\}\times[-b,b]$ is non empty and the
same holds for $\partial\Ome\cap\{a\}\times[-b,b]$,
\item $\partial\Ome\cap[-a,a]\times\{-b\}=\emptyset$ and $\partial\Ome
  \cap[-a,a] \times\{-b\}=\emptyset$.  
\end{enumerate}

We note $\Lambda$ the set of the closures of connected components of
$\partial\Ome\cap\inter{R_{a,b}}$ where
$\inter{R_{a,b}}=(-a,a)\times(-b,b)$. Let $\gamma\in\Lambda$ be one of these
boundary components, $\gamma$ is homeomorphic either to a circle either to
$[0,1]$. If it is homeomorphic to a segment, either it joins
$\{-a\}\times[-b,b]$ to $\{a\}\times[-b,b]$ (by connectedness, there are
exactly two such components) or the two end points are on the same edge of
the rectangle $R_{a,b}$.

We need some more notations and remarks that we shall use in the following
proofs. Let 
$c:[0,1]\rightarrow R_{a,b}$ be a Lipschitz continuous path with
$c(0)\in[-a,a]\times\{-b\}$, $c(1)\in[-a,a]\times \{b\}$ and
$c(t)\in\inter{R_{a,b}}$ for $0<t<1$. We note $J_c$ the set of the connected
components of $[0,1]\cap c^{-1}(\overline{\Ome})$. Let $j\in J_c$, there exist
$e_j,o_j\in(0,1)$ such that $j=[e_j,o_j]$. There is a total order on
$J_c$. Let $j,j'\in J_c$, we note $j\lhd j'$ if $o_{j}<e_{j'}$; the order
$\unlhd$ is then a total order. We remark that $J_c$ has a minimum $j_{min}$
and a maximum $j_{max}$. We then have the following lemma.
\begin{lem}\label{lem}
Let $\Ome$ and $c$ be as above. Let $j\in J_c$ with $j\neq j_{min}$. We
consider $j'\lhd j$ and note $\gamma$ 
the element of $\Lambda$ to which $c(e_j)$ belongs. Then there exists $j''$
with $j'\unlhd j''\lhd j$ such that $c(o_{j''})$ belongs to $\gamma$.
\end{lem}
\begin{proof}
We note $j_0=\sup\{ i\in J_c\, |\, i\lhd j\}$. There is two
possibilities. First, $j_0\lhd j$, in this case $j_0\unrhd j'$ and
$c(o_{j_0},e_j)$ is a curve outside $\overline{\Ome}$. Because of the
different cases for $\gamma$, $c(o_{j_0})$ is then in $\gamma$. The second
possibility is $j_0=j$. This implies that there exists $i\in J_c$ with
$o_i<e_j$ and $e_j-o_i$ as small as we want. The point $c(e_j)$ is at a non
zero distance from the complementary of $\gamma$ in $\partial\Ome$. Since
$c$ is Lipschitz continuous, there exists $j'\unlhd i\lhd j$ with
$c(o_i)\in\gamma$.
\end{proof}

If $c$ is injective, the $c[e_j,o_j]$ are the connected components of
$c\big([0,1]\big)\cap\overline{\Ome}$. 

We note $\Delta_1$ the connected component of $R_{a,b}\backslash\Ome$ that
contains $(0,-b)$ and $\Delta_2$ the one that contains $(0,b)$. For
$i\in\{1,2\}$, we note $\gamma_i$ the element of $\Lambda$ that is included in
the boundary of $\Delta_i$. $\gamma_1$ and $\gamma_2$ are the two elements of
$\Lambda$ that are homeomorphic to a segment and join $\{-a\}\times[-b,b]$ to
$\{a\}\times[-b,b]$. 

We are then able to give our height estimate result.

\begin{thm}\label{estim}
Let $a>\frac{1}{H}$ and $b>0$ be real numbers. We consider $\Ome\in\R_{a,b}$
a domain with piecewise smooth boundary that satisfies the above hypotheses
$1.$, $2.$ and $3.$, we note $\Lambda$, $\gamma_1$ and $\gamma_2$ as
above. Let $\Lambda=\Lambda_1\cup\Lambda_2$ be a partition of $\Lambda$ such
that $\gamma_1\in\Lambda_1$ and $\gamma_2\in\Lambda_2$. For $i\in\{1,2\}$, we
note $\Gamma_i$ the part of the boundary
$\bigcup_{\gamma\in\Lambda_i}\gamma$. Let $u$ be a solution of 
\eqref{cmc} on $\Ome$ which is continuous on $\overline{\Ome}$. We then have
the following upper bound :
$$
d\big(F(\Gamma_1),F(\Gamma_2)\big)\le\frac{2}{H}
$$
with $d$ the distance for compact sets of $\R^2$ and $F$ defined as at the
begining of the section.
\end{thm}

First we shall prove a weaker version of this result 

\begin{th1'}
Let $a>\frac{1}{H}$ and $b>0$ be real numbers. We consider $\Ome\in\R_{a,b}$,
$\Lambda$ and $\Lambda=\Lambda_1\cup\Lambda_2$ a partition as in Theorem
\ref{estim}. For $i\in\{1,2\}$, we note
$\Gamma_i=\bigcup_{\gamma\in\Lambda_i}\gamma$. Let $u$ be a solution of
\eqref{cmc} on $\Ome$ which is continuous on $\overline{\Ome}$. We then have 
the following upper bound :
$$
d\big(F(\Gamma_1),F(\Gamma_2)\big)\le 2a
$$
\end{th1'}

\begin{proof}
The idea of the proof is that, if the estimate on the distance does not hold,
we would be able to get a sphere of radius $\frac{1}{H}$  
through the graph of $u$ and this is impossible by maximum principle. So let
us assume that the distance $d\big(F(\Gamma_1),F(\Gamma_2)\big)$ is greater
than $2a$.

The first part of the proof consist in finding the place where the sphere will
be located. 

Since $\gamma_1$ and $\gamma_2$ join $\{-a\}\times [-b,b]$ to $\{a\}\times
[-b,b]$ in $R_{a,b}$, $F(\gamma_1)$ and $F(\gamma_2)$ join $\{-a\}\times \R$
to $\{a\}\times \R$ in $[-a,a]\times \R$. Let $\gamma$ be in $\Lambda_1$ and
$(x,z)$ be a point in $F(\gamma)$. Since $d\big(F(\Gamma_1),
F(\Gamma_2)\big)>2a$, no point of $F(\Gamma_2)$  has $z$ as second
coordinate, then, if $\gamma'\in\Lambda_2$, $F(\gamma')$ is either above
$F(\gamma)$ (\textit{i.e.} $\min_{F(\gamma')}z\ge \max_{F(\gamma)}z$) or below
(\textit{i.e.} $\max_{F(\gamma')}z\le \min_{F(\gamma)}z$). Then $\gamma$
defines a partition $\Lambda_2=\Lambda_2^-(\gamma) \cup \Lambda_2^+(\gamma)$
with $\Lambda_2^-(\gamma)$ (resp. $\Lambda_2^+(\gamma)$) is the set of
$\gamma'\in\Lambda_2$ such that $F(\gamma')$ is below (resp. above)
$F(\gamma)$. In the same way, $\gamma\in\Lambda_2$ defines a partition
$\Lambda_1=\Lambda_1^-(\gamma) \cup \Lambda_1^+(\gamma)$.

In the following, we assume that $\gamma_1\in\Lambda_1^-(\gamma_2)$
($F(\gamma_1)$ is below $F(\gamma_2)$). If $\gamma_1\in\Lambda_1^+(\gamma_2)$,
the proof is the same by exchanging the labels $1$ and $2$.

We then define:
$$
u_1=\max\left\{u(x,y)\,\Big|\,(x,y)\in
  \bigcup_{\gamma\in\Lambda_1^-(\gamma_2)} \gamma  ,-\frac{1}{H} \le x
  \le\frac{1}{H}\right\} 
$$
We note $(x_1,y_1)\in\bigcup_{\gamma\in\Lambda_1^-(\gamma_2)} \gamma$ a point
such that $u(x_1,y_1)=u_1$ and note $g_1\in \Lambda_1^-(\gamma_2)$ the
boundary component that contains $(x_1,y_1)$. We then note:
$$
u_2=\min\left\{u(x_1,y)\,\Big|\,(x_1,y)\in
  \bigcup_{\gamma\in\Lambda_2^+(g_1)} \gamma\right\} 
$$
We remark that $u_2$ is well defined because $\gamma_2\in \Lambda_2^+(g_1)$
and, since $\gamma_2$ join $\{-a\}\times [-b,b]$ to $\{a\}\times [-b,b]$,
there is a point in $\gamma_2$ with first coordinate $x_1$. We have $u_2>u_1$
and: 

\begin{fact}{f1}
For all $z\in(u_1,u_2)$, there exists $y$ such that $(x_1,y)\in\Ome$ and
$u(x_1,y)=z$. 
\end{fact}

\medskip
\noindent Let us prove this fact. We consider $c:[-b,b]\rightarrow R_{a,b}$
defined by $c(t)=(x_1,t)$. We consider the set $J_c$ with its order. Let
$j_0\in J_c$ be such that $c(e_{j_0})$ or $c(o_{j_0})$ is $(x_1,y_1)\in
g_1$. We then note: 
$$
j_1=\min\{j\rhd j_0\,|\,u(c(o_j))\ge u_2\}
$$
The segment $j_1$ exists because $u(c(o_{j_{max}}))\ge u_2$ since
$c(o_{j_{max}})\in\gamma_2$. Besides $j_0\lhd j_1$. First let us prove that
$u(c(e_{j_1})\le u_1$. We note $\gamma$ the element of $\Lambda$ to which
$c(e_{j_1})$ belongs. By Lemma \ref{lem}, there exists $i\in J_c$ with
$j_0\unlhd i  
\lhd j_1$ such that $c(o_i)\in\gamma$. We have $u(c(o_i))<u_2$ by definition
of $j_1$. If $\gamma\in\Lambda_1$, $u(c(o_i))<u_2$ implies that $\gamma\in
\Lambda_1^-(\gamma_2)$ and then $u(c(e_{j_1}))\le u_1$ (definition of
$u_1$). If $\gamma\in\Lambda_2$, 
$u(c(o_i))<u_2$ implies that $\gamma$ belongs to $\Lambda_2^-(g_1)$
(definition of $u_2$) and $u(c(e_{j_1}))\le u_1$. Now, since
$c([e_{j_1},o_{j_1}])$ is connected and included in $\overline{\Ome}$,
$(u_1,u_2)\subset u\circ c(e_{j_1},o_{j_1})$ and this prove Fact \ref{f1}.

Let $t$ be in $\R$ and $D_t$ be the closed disk in $[-a,a]\times \R$ with
center $(0,u_1+t)$ and radius $\frac{1}{H}$. $D_0$ contains the point
$(x_1,u_1)$. The diameter of $D_t$ is $\frac{2}{H}$ which is less than $2a$
then we have: 
\begin{gather*}
D_t\cap F(\Gamma_1)\neq \emptyset\Longrightarrow D_t\cap F(\Gamma_2)=
\emptyset\\ 
D_t\cap F(\Gamma_2)\neq \emptyset\Longrightarrow D_t\cap F(\Gamma_1)=
\emptyset 
\end{gather*}
We define:
$$
t_0=\inf\left\{t>0\,|\,D_t\cap F(\Gamma_1)=\emptyset\right\}
$$
By compactness, $D_{t_0}\cap F(\Gamma_1)\neq \emptyset$ and then $D_t\cap
F(\Gamma_2)=\emptyset$ for $0\le t\le t_0$.

\begin{fact}{f2}
We have $u_1+t_0<u_2$.
\end{fact}

\medskip
\noindent Actually, if $u_1+t_0\ge u_2$, $t'=u_2-u_1$ is less than $t_0$ and
$D_{t'}$ contains the point of $F(\Gamma_2)$ that realizes $u_2$. This implies
that $D_{t'}\cap F(\Gamma_1)=\emptyset$ and contradicts the definition of
$t_0$. 

By compactness, there exists then $t_1>t_0$ such that $u_1+t_1<u_2$,
$D_{t_1}\cap F(\Gamma_1)= \emptyset$ and $D_{t}\cap F(\Gamma_2)=
\emptyset$ for all $0\le t\le t_1$.

\begin{fact}{f3}
Let $\gamma\in\Lambda$ be a boundary component, then there are no
$Z_1,Z_2\in\R$ 
such that there exist $X_1,X_2\in[-\frac{1}{H},\frac{1}{H}]$ with
$(X_i,Z_i)\in F(\gamma)$ and:
$$
Z_1<u_1+t_1-\sqrt{\frac{1}{H^2}-X_1^2}\le u_1+t_1+ \sqrt{\frac{1}{H^2}-X_2^2}
<Z_2 
$$
($F(\gamma)$ can not have points above and below the disk $D_{t_1}$)
\end{fact}

\medskip
\noindent First we suppose that $\gamma\in\Lambda_1$. Since $Z_2>u_1$, the
definition of $u_1$ implies that $\gamma$ belongs to
$\Lambda_1^+(\gamma_2)$. Since $\gamma_2$ joins $\{-a\}\times 
[-b,b]$ to $\{a\}\times[-b,b]$, $F(\gamma_2)$ has a point of coordinates
$(x_1,z)$; by definition of $u_2$, $z\ge u_2$. Then the second coordinate of
every point of $F(\gamma)$ needs to be more than $u_2$ this contradicts
$Z_1<u_1+t_1$. Now if $\gamma\in\Lambda_2$, $D_t$ does not intersect
$F(\Gamma_2)$ for $0\le t\le t_1$, so letting go down the disk
$D_t$ from $t_1$ to $0$, we get $Z_1\le u_1$ and then $\gamma$ belongs to
$\Lambda_2^-(g_1)$. This implies that the second coordinate of
every point of $F(\gamma)$ needs to be less than $u_1$ this contradicts
$Z_2>u_1+t_1$ and proofs Fact \ref{f3}. 

The idea is now to get a sphere of radius $\frac{1}{H}$ through the
disk $D_{t_1}$. We note $S_v$ the sphere of radius $\frac{1}{H}$ and center
$(0,v,u_1+t_1)$. When $v$ changes, $S_v$ moves in an horizontal cylinder with
vertical section $D_{t_1}$. For far from zero negative $v$, $S_v$ is out
$\overline{\Ome}\times\R$. Since $D_{t_1}\cap F(\Gamma_1)= \emptyset$ and
$D_{t_1}\cap 
F(\Gamma_2)= \emptyset$, $S_v$ does not intersect the boundary of the graph of
$u$ for any $v$. The graph of $u$ splits $\overline{\Ome}\times\R$ into two
connected components: $\boG^+$, above the graph, and $\boG^-$ which is
below. Since $u_1\le u_1+t_1\le u_2$, there exists $v$ such that $S_v$
intersects the  graph of $u$ from Fact \ref{f1}.

We start with far from zero negative $v$ and let $v$ increase until $v_0$
which is the first contact between the graph and the sphere. This first
contact does not occur in the boundary of the graph since the sphere never
intersects it. Then, since the graph is not a piece of a sphere  because of
the size of $\Ome$, the maximum principle implies that, in the neighborhood of
the contact point, the sphere $S_{v_0}$ is in $\boG^-$ (we recall that the
mean curvature vector of the graph points in $\boG^+$ because of the equation
\eqref{cmc}). But, in fact, we have:

\begin{fact}{f4}
In the neighborhood of the contact point, the sphere $S_{v_0}$ is in $\boG^+$.
\end{fact}

\medskip
\noindent This fact is not clear since, because of the shape of the domain
$\Ome$, the sphere do not stop to get in and out $\overline{\Ome}\times \R$. We
note $p=(x,y,z)$ the first contact point; we know that 
$(x,z)\in D_{t_1}$. We define $c:s\mapsto (x,s)\in R_{a,b}$ and consider
$J_c$. We have $(x,s,z)\in\boG^-$ for $s<y$ near $y$ since the sphere
$S_{v_0}$ is in $\boG^-$ in a neighborhood of $p$. Then there exists:
$$
s_0=\min\left\{s\le -b\,|\, (x,s,z)\in\boG^-\right\}
$$
Since $p$ is the first contact point, there is $j\in J_c$ such that
$s_0=e_j$. We know that $D_{t_1}$ is above $F(\gamma_1)$ then
$(x,e_{j_{min}},z)\in\boG^+$, so $j\rhd j_{min}$. We know then that there
exists $j'\lhd j$ such that $c(e_j)$ and $c(o_{j'})$ belong to
the same element of $\Lambda$; then, if $(x,e_j,z)\in\boG^-$,
$(x,o_{j'},z)\in\boG^-$, this is due to Fact \ref{f3}. This implies that
$(x,s,z)\in\boG^-$ for $s<o_{j'}$ near $o_{j'}$ and then $s_0<e_j$; we have
our contradiction. 

\medskip
This end the proof of $d\big(F(\Gamma_1),F(\Gamma_2)\big)\le 2a$
\end{proof}

We remark that Theorem 2' will be sufficient for most of the applications and
Theorem \ref{estim} is just an improvement. So let us replace $2a$ by
$\frac{2}{H}$ to get Theorem \ref{estim}. 

\begin{proof}[Proof of Theorem \ref{estim}]
Let us consider $a'>0$ with $\frac{1}{H}<a'<a$. The idea of the proof is to
apply Theorem 2' to a well chosen connected component of
$\Ome\cap R_{a',b}$. We note $D^i$ the connected components of
$\Ome\cap R_{a',b}$. First we remark that, among these components, there are
ones that satisfy the hypotheses $1.$, $2.$ and $3.$; for example, since 
$\gamma_1$ joins $\{-a\}\times[-b,b]$ to $\{a\}\times[-b,b]$, one connected
component of $\gamma_1\cap R_{a',b}$ joins $\{-a'\}\times[-b,b]$ to
$\{a'\}\times[-b,b]$ then a $D^i$ that has this component in its boundary
satisfies the three hypotheses. A component of $\Ome\cap R_{a',b}$ that
satisfies the hypotheses is called a good component and the other ones are the
bad ones; we rename these good components $D^1,\dots,D^k$. There is only a
finite number of such components since the length of the part of $\partial
\Ome$ in $R_{a',b}$ is finite.

Let us consider a good component $D^i$. As defined at the beginning of the
section, a set of boundary component $\Lambda^i$ is associated to $D^i$. In
$\Lambda^i$ there is two particular elements, these are the two boundary
components which are homeomorphic to a segment and joins $\{-a'\}\times[-b,b]$
to $\{a'\}\times[-b,b]$. To avoid any confusion, we note these components
$\gamma_\alpha^i$ and $\gamma_\beta^i$ ($\gamma_\alpha^i$ is a part of the
boundary of the connected component of $R_{a',b}\backslash D^i$ that contains
$(0,-b)$ and $\gamma_\beta^i$ is the other one). Each element of $\Lambda^i$
is a part of an element of $\Lambda$, then we get a partition
$\Lambda^i=\Lambda_1^i\cup \Lambda_2^i$: an element of $\Lambda^i$ is in
$\Lambda_1^i$ (resp. $\Lambda_2^i$) if it is a part of a element of
$\Lambda_1$ (resp. $\Lambda_2$). Now the proof consists in applying Theorem
\ref{estim}' to a component $D^i$ such that $\gamma_\alpha^i\in\Lambda_1^i$
and  $\gamma_\beta^i \in \Lambda_2^i$.

To each good component $D^j$, we can associate a real number which is the
second coordinate of the end point of $\gamma_\alpha^j$ in $\{-a\}\times
[-b,b]$. In the folowing, we order the good components with respect to this
real number and rename the good components $D^1,\dots,D^k$ with respect to
this order. The order is the same if we consider the second coordinate of 
$\gamma_\alpha^j \cup \{a\}\times [-b,b]$, $\gamma_\beta^j \cup \{-a\}\times
[-b,b]$ or $\gamma_\beta^j \cup \{a\}\times [-b,b]$.

Let $D$ be a bad component of $\Ome\cap R_{a',b}$. In fact, it is a bad
component only because of the hypothesis $2.$, then, in $D$, there is no path
from $\{-a'\}\times [-b,b]$ to $\{a'\}\times [-b,b]$. This implies that, as
in Figure \ref{courbec}, there
exists a path $c:[0,1]\rightarrow R_{a',b}$ that joins $(0,-b)$ to $(0,b)$,
is outside all the bad components and such that there exist $0<e_1<o_1<e_2<
\cdots<e_k<o_{k}<1$ with: 
$$
c([0,1])\cap\Ome=\bigcup_i c(e_i,o_i)
$$
and $c(e_i,o_i)\subset D^i$. First we remark that $c(e_1)$ is in
$\gamma_1$ so it is in  $\gamma_\alpha^1\in\Lambda_1^1$ and $c(o_k)$
is in $\gamma_2$ so it is in $\gamma_\beta^k\in\Lambda_2^k$. Then there
exists:  
$$
i_0=\min\left\{i\,\Big|\,\gamma_\beta^i \in \Lambda_2^i\right\} 
$$

\begin{figure}[h]
\begin{center}
\resizebox{0.6\linewidth}{!}{\input{figestimH2.pstex_t}}
\caption{\label{courbec}}
\end{center}
\end{figure}

Let us prove that the good component $D^{i_0}$ is such that
$\gamma_\alpha^{i_0}\in \Lambda_1^{i_0}$. First we assume that it is not the
case, \textit{i.e.} $\gamma_\alpha^{i_0}\in \Lambda_2^{i_0}$. Since
$\gamma_\alpha^{1}\in \Lambda_1^{1}$, $i_0>0$. $\gamma_\alpha^{i_0}$ is part
of an element $\gamma$ of $\Lambda_2$ then $c(e_{i_0})\in\gamma_\alpha^{i_0}
\subset \gamma$. Besides we know, by Lemma \ref{lem}, that $c(o_{i_0-1})$
belongs to the same component as $c(e_{i_0})$ then $c(o_{i_0-1})\in \gamma$
and $\gamma_\beta^{i_0-1}\subset \gamma$; this implies that
$\gamma_\beta^{i_0-1} \in \Lambda_2^{i_0-1}$. This is a contradiction with the
definition of $i_0$.  

\medskip
Now, we apply Theorem 2' to $D^{i_0}$. We note, for $j\in\{1,2\}$,
$\Gamma_j'=\bigcup_{\gamma\in\Lambda_j^{i_0}}\gamma$. Then we have
$\gamma_\alpha^{i_0}\in\Gamma_1'$ and $\gamma_\beta^{i_0}\in\Gamma_2'$; so we
can apply the theorem and we get:
$$
d\big(F(\Gamma_1'),F(\Gamma_2')\big)\le 2a'
$$
Besides, we have $\Gamma_1'\subset\Gamma_1$ and $\Gamma_2'\subset \Gamma_2$
then:
$$
d\big(F(\Gamma_1),F(\Gamma_2)\big)\le d\big(F(\Gamma_1'),F(\Gamma_2')\big)\le
2a' 
$$

This inequality is true for every $a'>\frac{1}{H}$, so:
$$
d\big(F(\Gamma_1),F(\Gamma_2)\big) \le \frac{2}{H}
$$
\end{proof}


\section{Some consequences of Theorem \ref{estim}}
\label{conseq}
The aim of this section is to give some consequences of Theorem \ref{estim}
for solutions of the Dirichlet problem associated to the constant mean
curvature equation \eqref{cmc} on unbounded domains.

First we explain what kind of domains we shall consider. Let $b_-$ and
$b_+:\R_+\rightarrow \R$ be two continuous functions such 
that, for every $x\ge 0$, $b_-(x)<b_+(x)$. We are interested in domains of the
type $\Ome=\{(x,y)\in\R_+\times \R\,|\, b_-(x)<y<b_+(x)\}$. When a solution
$u$ is solution of \eqref{cmc} on $\Ome$ and continuous on $\overline{\Ome}$,
we define two continous functions $f_-$ and $f_+$ by
$f_\pm(x)=u(x,b_\pm(x))$. $f_-$ and $f_+$ are the boundary values of $u$.

Let us fix a last definition, if $x\in\R_+$, we note
$I_x=\{x\}\times[b_-(x),b_+(x)]$. We then have the following height estimate:
\begin{prop}\label{minor}
Let $b_-$and $b_+$ be continuous functions on $\R^+$ with $b_-(x)<b_+(x)$. We
note $\Ome=\{(x,y)\in\R_+\times \R\,|\, b_-(x)<y<b_+(x)\}$. We consider $u$ a
solution of \eqref{cmc} on $\Ome$ continuous on $\overline{\Ome}$. We consider
$x_0>\frac{2}{H}$ and $M$ such that: 
\begin{equation}\label{equati}
\min_{[x_0-\frac{2}{H},x_0+\frac{2}{H}]} (f_-,f_+)\ge M 
\end{equation}
Then:
$$
\min_{I_{x_0}}u\ge M-\frac{3}{H}
$$
\end{prop}

\begin{proof}
Let $\eps$ be a positive number. Since $f_-$ and $f_+$ are
continuous, there exist $\eta>0$ such that:
$$
\min_{[x_0-\frac{2}{H}-\eta,x_0+\frac{2}{H}+\eta]} (f_-,f_+)\ge M-\eps
$$
Let us now suppose that there is a piecewise smooth injective path $c:
[0,1]\rightarrow 
\Ome\bigcap[x_0-\frac{2}{H}-\eta,x_0] \times\R$ that joins 
$I_{x_0-\frac{2}{H}-\eta}$ to $I_{x_0}$ such that:
\begin{equation}\label{eq3}
u\circ c(t)<M-\eps-\frac{2}{H}
\end{equation}

We consider the domain $D$ bounded by $c$, the curve $y=b_-(x)$ for $x\in
[x_0-\frac{2}{H}-\eta,x_0]$, a segment included in $I_{x_0-\frac{2}{H}-\eta}$
and one included in $I_{x_0}$; $D$ statifies the three hypotheses of Section
\ref{estimH1}. In fact, since the function $b_-$ is only continuous the
boundary of $D$ is not piecewise smooth, but Theorem \ref{estim} can be
applied because of the shape of $D$. Then because of \eqref{equati} and
\eqref{eq3},  Theorem \ref{estim} is not satisfied. Finally, this implies that
there exists a path $c_1: [0,1]\rightarrow
\Ome\bigcap[x_0-\frac{2}{H}-\eta,x_0] \times\R$ that joins the curve
$y=b_-(x)$ to the curve $y=b_+(x)$ such that $u\circ c_1(t)\ge
M-\eps-\frac{2}{H}$. 

By the same arguments, there exists a path $c_2:[0,1]\rightarrow
\Ome\bigcap[x_0,x_0+\frac{2}{H}+\eta] \times\R$ that joins the curve
$y=b_-(x)$ to the curve $y=b_+(x)$ such that $u\circ c_2(t)\ge
M-\eps-\frac{2}{H}$. 

\begin{figure}[h]
\begin{center}
\resizebox{0.8\linewidth}{!}{\input{figestimH1.pstex_t}}
\caption{\label{domaineD}}
\end{center}
\end{figure}

Now the domain $D$ bounded by $c_1$, $c_2$, a piece of $y=b_-(x)$ and a piece
of $y=b_+(x)$ contains $I_{x_0}$ (see Figure \ref{domaineD}). Besides on the
boundary of $D$, $u$ is 
everywhere greater than $M-\eps-\frac{2}{H}$ by \eqref{equati} and above. Then
by a classical height estimate \cite{Se}, $u$ is greater than
$M-\eps-\frac{3}{H}$ in $D$. this gives us:  
\begin{equation}\label{eq4}
\min_{I_{x_0}}u\ge M-\eps-\frac{3}{H}
\end{equation}
Letting $\eps$ goes to zero, we get the expected result.
\end{proof}

We have also a simple upper-bound in this case.
\begin{prop}\label{major}
Let $b_-$ and $b_+$ be continuous functions on $\R^+$ with $b_-(x)<b_+(x)$. We
note $\Ome=\{(x,y)\in\R_+\times \R\,|\, b_-(x)<y<b_+(x)\}$. We consider $u$ a
solution of \eqref{cmc} on $\Ome$ continuous on $\overline{\Ome}$. We consider
$x_0>\frac{1}{2H}$ and $M$ such that: 
$$
\max_{[x_0-\frac{1}{2H},x_0+\frac{1}{2H}]} (f_-,f_+)\le M 
$$
Then:
$$
\max_{I_{x_0}}u\le M
$$
\end{prop}

\begin{proof}
Let $\eps$ be a positive number. Since $f_-$ and $f_+$ are continuous, there
exist $\eta>0$ such that: 
$$
\max_{[x_0-\frac{1}{2H}-\eta,x_0+\frac{1}{2H}+\eta]} (f_-,f_+)\le M+\eps
$$
Let us consider the cylinder of radius $\frac{1}{2H}$ which is centered on
the axis $\{x=x_0\}\cap \{z=t\}$.  For big $t$, the cylinder is above the graph
and we can let $t$ decrease. Until $t=M+\frac{1}{2H}$, the cylinders can not
meet the boundary of the graph of $u$. The maximum principle then says us
that, for $t=M+\frac{1}{2H}$, the cylinder is still above the graph; so we
get: 
$$
\max_{I_{x_0}}u\le M+\epsilon
$$
Letting $\eps$ goes to zero, we get the expected result.
\end{proof}

Let $f: I\rightarrow \R$ be a function, we define the variation of $f$
around the point $x_0$ by:
$$
V_t(x_0,f)=\sup_{[x_0-t,x_0+t]}f -\inf_{[x_0-t,x_0+t]}f
$$
Let $f$ and $g$ be two continuous functions $I\rightarrow \R$; we
define the variation of the pair $(f,g)$ around $x_0$ by:
$$
V_t(x_0,f,g)=\max\left(V_t(x_0,f), V_t(x_0,g)\right)
$$

The two preceding propositions give us the following result:
\begin{thm}\label{th3}
Let $b_-$ and $b_+$ be continuous functions on $\R^+$ with $b_-(x)<b_+(x)$. We
note $\Ome=\{(x,y)\in\R_+\times \R\,|\, b_-(x)<y<b_+(x)\}$. We consider $u$ a
solution of \eqref{cmc} on $\Ome$ continuous on $\overline{\Ome}$. We consider
$x_0>\frac{2}{H}$ and $M$ such that: 
$$
V_{\frac{2}{H}}(x_0,f_-,f_+)\le M
$$
Then there exists $M'$ which depends only on $M$ and $H$ such that:
$$
\max_{I_{x_0}}u-\min_{I_{x_0}}u\le M'
$$
For example, $M'=4M+\frac{5}{H}$ works.
\end{thm}

\begin{proof}
We have for $x\in[x_0-\frac{2}{H},x_0+\frac{2}{H}]$,
$|f_-(x)-f_-(x_0)|\le M$ and $|f_+(x)-f_+(x_0)|\le M$. Then if we
apply Theorem \ref{estim} to $\Ome\bigcap [x_0-\frac{2}{H},
x_0+\frac{2}{H}]\times \R$, we get that the graph of $f_-$ over this
segment is at a distance less than $\frac{2}{H}$ from the one of $f_+$. Since,
for $\alpha\in\{-,+\}$, the graph of $f_\alpha$ is in the horizontal strip
$f_\alpha(x_0)-M\le z\le f_\alpha(x_0)+M$, we have 
$|f_-(x_0)-f_+(x_0)|\le 2(M+\frac{1}{H})$. This implies that, for every
$x,x'\in [x_0-\frac{2}{H},x_0+\frac{2}{H}]$, $|f_\alpha(x)-f_\beta(x')|\le
4M+\frac{2}{H}$ with $\alpha,\beta\in\{-,+\}$. Then there exist $A\in\R$ such
that, for every $x\in[x_0-\frac{2}{H},x_0+\frac{2}{H}]$, we have:
\begin{gather}
\label{eq1}
|f_-(x)-A|\le 2M+\frac{1}{H}\\
\label{eq2}
|f_+(x)-A|\le 2M+\frac{1}{H}
\end{gather}
These two equations with Proposition \ref{minor} implies that:
\begin{equation}\label{minoration}
\min_{I_{x_0}}u\ge A-2M-\frac{4}{H}
\end{equation}
With Proposition \ref{major}, we get:
\begin{equation}\label{majoration}
\max_{I_{x_0}}u\le A+2M+\frac{1}{H}
\end{equation}
Then, in bringing together \eqref{minoration} and \eqref{majoration}, we
obtain: 
$$
\max_{I_{x_0}}u-\min_{I_{x_0}}u\le 4M+\frac{5}{H}
$$
\end{proof}

Theorem  \ref{th3} has an easy corollary.
\begin{cor}\label{co}
Let $b_-$ and $b_+$ be continuous functions on $\R^+$ with $b_-(x)<b_+(x)$. We 
note $\Ome=\{(x,y)\in\R_+\times \R\,|\, b_-(x)<y<b_+(x)\}$. We consider $u$ a
solution of \eqref{cmc} on $\Ome$ continuous on $\overline{\Ome}$. We consider
$x_0>\frac{2}{H}$ and $M$ such that: 
$$
V_{\frac{2}{H}}(x_0,f_-,f_+)\le M
$$

Then there exists $M'$ which depends only on $M$ and $H$ such that, for every
$p\in I_{x_0}$ and $\alpha\in\{-,+\}$, we have:
$$
f_\alpha(x_0)- M'\le u(p)\le f_\alpha(x_0)+ M'
$$
For example, $M'=4M+\frac{5}{H}$ works. 
\end{cor}
\begin{proof}
It is just the fact that $(x_0,b_-(x_0))$ and $(x_0,b_+(x_0))$ are in
$I_{x_0}$. 
\end{proof}


\section{Two uniqueness results}
In this section, we use Corollary \ref{co} to prove uniqueness theorems for
the Dirichlet problem associated to \eqref{cmc}.

\begin{thm}\label{uni1}
Let $b_-,b_+$ be two continuous functions on $\R_+$ such that
$b_-(0)=b_+(0)$ and $b_-(x)<b_+(x)$ for every $x>0$. We note
$\Ome=\{(x,y)\in\R_+\times\R\, |\, b_-(x)<y<b_+(x)\}$. Let $f_-,f_+$ be two
continuous functions on $\R_+$ such that $f_-(0)=f_+(0)$. We suppose
that there exist an increasing positive sequence $(x_n)_{n\in\N}$ with
$\lim x_n= +\infty$ and a sequence $(M_n)_{n\in\N }$ with $M_n=o(\ln x_n)$
such that, for every $n\in\N$, we have:
$$
V_{\frac{2}{H}}(x_n,f_-,f_+)\le M_n
$$
Then, if there exists a solution $u$ of \eqref{cmc} on $\Ome$ with value $f_-$
and $f_+$ on the boundary, this solution is unique. 
\end{thm}

\begin{proof}
Let $u_1$ and $u_2$ be two different solutions of the Dirichlet problem with
$f_-$ and $f_+$ as boundary data. We know by a result of P.~Collin and
R.~Krust \cite{CK} that: 
$$
\liminf_{x\rightarrow+\infty}\frac{\max_{I_x}|u_1-u_2|}{\ln x}>0
$$ 

But by corollary \ref{co}, we know that:
\begin{gather*}
\max_{I_{x_n}}|u_1-f_-(x_n)|\le 4M_n+\frac{5}{H}\\
\max_{I_{x_n}}|u_2-f_-(x_n)|\le 4M_n+\frac{5}{H}
\end{gather*}
So, we get:
$$
\max_{I_{x_n}}|u_1-u_2|\le 8M_n+\frac{10}{H}
$$
By the hypothesis on $M_n$, we have:
$$
\lim_{n\rightarrow\infty}\frac{\max_{I_{x_n}}|u_1-u_2|}{\ln x_n}=0
$$
This gives us a contradiction since $x_n\rightarrow+\infty$.
\end{proof}

We have also a second theorem.
\begin{thm}\label{uni2}
Let $b_-,b_+$ be two continuous functions on $\R$ such that
$b_-(x)<b_+(x)$ for every $x\in\R$. We note
$\Ome=\{(x,y)\in\R^2\, |\, b_-(x)<y<b_+(x)\}$. Let $f_-,f_+$ be two
continuous functions on $\R$. We suppose
that there exist one increasing sequence $(x_n)_{n\in\N}$
with $\lim x_n= +\infty$ and one decreasing sequence $(x_n')_{n\in\N}$ with
$\lim x_n'=-\infty$ and two sequences $(M_n)_{n\in\N}$ and
$(M_n')_{in\in\N}$ such that $M_n=o(\ln |x_n|)$, $M_n'=o(\ln |x_n'|)$ and, for
every $n\in\N$, we have: 
\begin{gather*}
V_{\frac{2}{H}}(x_n,f_-,f_+)\le M_n\\
V_{\frac{2}{H}}(x_n',f_-,f_+)\le M_n'
\end{gather*}
Then, if there exists a solution $u$ of \eqref{cmc} on $\Ome$ with value $f_-$
and $f_+$ on the boundary, this solution is unique. 
\end{thm}

\begin{proof}
Let $u_1$ and $u_2$ be two different solutions of the Dirichlet problem with
$f_-$ and $f_+$ as boundary data. We know (\cite{CK}) that:
\begin{equation}\label{eqE}
\max_{I_x\cup I_{-x}}|u_1-u_2|\xrightarrow[x\rightarrow +\infty]{} +\infty 
\end{equation}
We note $C=\max_{I_0}|u_1-u_2|$. We then have $u_2-C-1< u_1<u_2+C+1$ on $I_0$
and because of \eqref{eqE} the set $\{|u_1-u_2|>M+1\}$ is non-empty. then we
can assume that there exists a subdomain $\Ome^*\subset\Ome\cap\R^-\times\R$
which is a connected component of $\{u_1>u_2+M+1\}$. Then we have:
$$
\liminf_{x\rightarrow-\infty}\frac{\max_{I_x\cap\Ome^*}|u_1-u_2-C-1|} {\ln
  |x|}>0 
$$

As in the preceding proof, Corollary \ref{co} says us that, for every $n$, we
have: 
$$
\max_{I_{x_n'}}|u_1-u_2-C-1|\le 8M_n'+\frac{10}{H}+C+1
$$
By the hypothesis on $M_n'$, we have:
$$
\lim_{n\rightarrow\infty}\frac{\max_{I_{x_n'}}|u_1-u_2-C-1|}{\ln |x_n'|}=0
$$
This gives us a contradiction since $x_n'\rightarrow-\infty$ and ends the
proof. 
\end{proof}

This theorem can be used to study the uniqueness of the solutions which were
built by P.~Collin in \cite{Co} and by R.~Lop\'ez in \cite{Lo}.

There are others results of uniqueness we can prove with the same
arguments. For example, if we suppose that $b_+-b_-$ is bounded in Theorems
\ref{uni1} and \ref{uni2}, we need only to assume that $M_n=o(x_n)$ and
$M_n'=o(|x_n'|)$ to have the uniqueness.


\bigskip

\noindent Laurent Mazet

\noindent Laboratoire Emile Picard (UMR 5580), Universit\'e Paul Sabatier,

\noindent 118, Route de Narbonne, 31062 Toulouse, France.

\noindent E-mail: mazet@picard.ups-tlse.fr

\end{document}